\documentclass[a4paper,11pt,twoside]{article}

\RequirePackage[a4paper,head=1cm]{geometry}
\usepackage[utf8x]{inputenc}
\usepackage{amsmath}
\usepackage{amssymb}
\usepackage{hyperref}
\usepackage{dsfont}
\newtheorem{Definition}{Definition}[section]
\newtheorem{Theoreme}{Theorem}
\newtheorem{Lemme}{Lemma}[section]
\newtheorem{Corollaire}{Corollary}[section]
\newtheorem{Proposition}{Proposition}[section]
\newtheorem{Remarque}{\bf Remark}

\title{Improved Sobolev Inequalities and Muckenhoupt weights \\ on stratified Lie groups } 
\author{Diego Chamorro } 
\begin{document} 

\maketitle 
\begin{scriptsize}
\abstract{We study in this article the Improved Sobolev inequalities with Muckenhoupt weights within the framework of stratified Lie
groups. This family of inequalities estimate the $L^{q}$ norm of a function by the geometric mean of two norms corresponding to Sobolev spaces $\dot{W}^{s,p}$ and Besov spaces $\dot{B}^{-\beta, \infty}_\infty$. When the value $p$ which characterizes Sobolev space is strictly larger than $1$, the required result is well known in $\mathbb{R}^n$ and is classically obtained  by a Littlewood-Paley dyadic blocks manipulation. For these inequalities we will develop here another totally different technique. When $p=1$, these two techniques are not available anymore and following M. Ledoux in \cite{Ledoux} in $\mathbb{R}^n$, we will treat here the critical case $p=1$ for general stratified Lie groups in a weighted functional space setting. Finally, we will go a step further with a new generalization of Improved Sobolev inequalities using weak-type Sobolev spaces.\\[3mm] \textbf{Keywords:} Improved Sobolev inequalities, stratified Lie groups, Muckenhoupt weights.}\end{scriptsize}


\section{Introduction}

In the Euclidean case, we can roughly distinguish three types of Improved Sobolev inequalities following the method used in their proof and the parameter's range defining the functional spaces. Let us recall these inequalities (for a precise definition of the functional spaces used below, please refer to section \ref{espa}).\\

Historically the first method, due to P. G\'erard, F. Oru and Y. Meyer \cite{GMO}, is based on a Littlewood-Paley decomposition and interpolation results applied to dyadic blocks. For a function $f$ such that $f\in \dot{W}^{s_1,p}(\mathbb{R}^n)$ and $f\in \dot{B}^{-\beta,\infty}_{\infty}(\mathbb{R}^n)$, the inequality obtained reads as follows:
\begin{equation}\label{ISI1}
\|f\|_{\dot{W}^{s,q}}\leq C \|f\|_{\dot{W}^{s_1,p}}^\theta\|f\|_{\dot{B}^{-\beta,\infty}_{\infty}}^{1-\theta}
\end{equation}
where $1<p<q<+\infty$, $\theta=p/q$, $s=\theta s_1 -(1-\theta)\beta$ and $-\beta<s<s_1$. Let us stress that the value $p=1$ is forbidden here. We write $\dot{W}^{s,p}$ for homogeneous $(s,p)$-Sobolev spaces and  $\dot{B}^{-\beta,\infty}_{\infty}$ for homogeneous $(-\beta,\infty, \infty)$-Besov spaces.\\

The second method, studied by M. Ledoux in \cite{Ledoux}, use semi-group properties related to Laplacian and heat kernel and allows us to treat the case $p=1$. If  $\nabla f\in L^{p}(\mathbb{R}^n)$ and $f\in \dot{B}^{-\beta,\infty}_{\infty}(\mathbb{R}^n)$, we have
\begin{equation}\label{ISI2}
\|f\|_{L^q}\leq C \|\nabla f\|_{L^p}^\theta\|f\|_{\dot{B}^{-\beta,\infty}_{\infty}}^{1-\theta}
\end{equation}
with $1\leq p < q<+\infty$, $\theta=p/q$ and $\beta=\theta/(1-\theta)$.\\ 

Finally, the third method proposed by A. Cohen, W. Dahmen, I. Daubechies \& R. De Vore in \cite{Cohen2} use a BV-norm weak estimation using wavelet coefficients and isoperimetric inequalities and gives, for a function $f$ such that $f\in BV(\mathbb{R}^n)$ and $f\in \dot{B}^{-\beta,\infty}_{\infty}(\mathbb{R}^n)$, the estimation below:
\begin{equation}\label{ISI3}
\|f\|_{\dot{W}^{s,q}}\leq C \| f\|_{BV}^{1/q}\|f\|_{\dot{B}^{-\beta,\infty}_{\infty}}^{1-1/q}
\end{equation}
where $1<q\leq 2$, $0\leq s<1/q$ and $\beta=(1-sq)/(q-1)$. When $s=0$, this last result implies (\ref{ISI2}) with $p=1$, but is limited by the fact that $1<q\leq 2$.\\

In this paper we will study these inequalities using as a framework stratified Lie groups which are a natural generalization of $\mathbb{R}^n$ when modifying dilations (see (\ref{DefDilat}) below for a definition). However, this seemingly simple modification induces some serious technical problems at many levels since the whole group structure is changed: for example, the underlying geometry is totally different (see \cite{Strichartz} and \cite{Garofalo}) which makes the techniques used in \cite{Cohen2} hardly transposable to this setting; observe also that the use of Fourier transform, and the classical associated tools, is not as straightforward as in $\mathbb{R}^n$.\\

For the Heisenberg group, inequalities of type (\ref{ISI1}) have been carried out in \cite{Bahouri} and with the work realized in \cite{Furioli2} we can deduce these inequalities for stratified Lie groups in a unweighted setting. Note that these authors develop systematically in each case a Littlewood-Paley decomposition in order to obtain these estimates, we will show here how to treat these inequalities in a far more direct and simpler way using maximal functions.\\ 

This is one of the main novelties of this paper, however, our principal aim is to generalize inequality (\ref{ISI2}) to stratified Lie groups using weighted spaces and to give a new weak-type estimation which lies, roughly speaking, between (\ref{ISI2}) and (\ref{ISI3}). In order to achieve this, we will develop some techniques using properties associated to the sub-Laplacian spectral decomposition. \\

We will consider throughout this paper weighted functional spaces with weights $\omega$ belonging to the Muckenhoupt classes $A_{p}$ for $1\leq p<+\infty$\footnote{When $p=+\infty$, we do not consider them, since in this case weighted spaces $X^{\infty}(\omega)$ coincide with traditional ones $X^{\infty}$ where $X$ is a Lebesgue, Sobolev or Besov space.}. The main reason for considering these weights lies in their connection with maximal functions which will lead us to a painless proof for inequalities of type (\ref{ISI1}).\\

Our principal theorem treats the critical case $p=1$ of improved Sobolev inequalities: 
\begin{Theoreme}\label{smile00} 
Let $\mathbb{G}$ be a stratified Lie group and $\omega$ a weight in the Muckenhoupt class $A_{1}$.  If $\nabla f\in L^{1}(\mathbb{G},\omega)$ and $f\in \dot{B}^{-\beta, \infty}_{\infty}(\mathbb{G})$, then we have the following inequalities:
\begin{enumerate}
\item[$\bullet$] \textbf{\emph{[Strong inequalities]}}
\begin{equation}\label{smile0} 
\|f\|_ {L^{q}(\omega)}\leq C \|\nabla f \|_ {L^{1}(\omega)}^{\theta} \|f\|_{\dot{B}^{-\beta, \infty}_{\infty}}^{1-\theta} 
\end{equation} 
where $1<q<+\infty$, $\theta = 1/q$ and $\beta=\theta/(1-\theta)$.  

\item[$\bullet$] \textbf{\emph{[Weak inequalities]}} 
\begin{equation}\label{smile1} 
\|f\|_ {\dot{W}^{s, q}_{\infty}(\omega)}\leq C \|\nabla f\|_{L^1(\omega)}^{\theta}\|f\|_{\dot{B}^{-\beta,\infty}_{\infty}}^{1-\theta} 
\end{equation} 
where $1<q<+\infty$, $0<s<1/q<1$, $\theta=1/q$ and $\beta=\frac{1-sq}{q-1}$.
Here $\|\cdot\|_{\dot{W}^{s, q}_{\infty}(\omega)}$ characterizes the weighted weak homogenous Sobolev space which definition is given by formula (\ref{DefDebilSoboLe}) in section \ref{espa}.
\end{enumerate}
\end{Theoreme} 

It is possible to see inequality (\ref{smile1}) as a weak-type improvement of (\ref{ISI3}). Indeed, the general Sobolev-like inequality obtained in \cite{Cohen2} uses in fact a Besov space in the left-hand side: 
$$\|f\|_{\dot{B}^{s,q}_q}\leq C \| f\|_{BV}^{1/q}\|f\|_{\dot{B}^{-\beta,\infty}_{\infty}}^{1-1/q}$$
which turns to be (\ref{ISI3}) only if $1<q\leq 2$ since in this case we have $\dot{B}^{s,q}_q\subset \dot{W}^{s,q}$. The restriction $q\in ]1,2]$ is a serious limitation as, for $q>2$, this Besov-Sobolev spaces embbeding is reversed. It is then important to observe that our weak inequality does not have this restriction since we allow $q$ to be in the interval $]1, +\infty[$. \\

Note also that in the context of stratified Lie groups, the weak inequality (\ref{smile1}) is the sharpest result available.\\

Our second result provides the main tool for proving theorem \ref{smile00}:
\begin{Theoreme}[Modified pseudo-inequality of Poincar\'e]\label{CHAME} 
Let $\mathbb{G}$ be a stratified Lie group, $\omega \in A_{1}$ and $\nabla f\in L^{1}(\mathbb{G}, \omega)$.  We have the following estimate for $0\leq s<1$ and for $t>0$:  
\begin{equation}\label{CHAME1} 
\|\mathcal{J}^{s/2}f-H_{t}\mathcal{J}^{s/2}f \|_ { L^{1}(\omega)}\leq C \; t^{\frac{1-s}{2 } } \|\nabla f\|_ { L^{1}(\omega) }.  
\end{equation}
Where $\mathcal{J}$ is a sub-Laplacian on $\mathbb{G}$ invariant with respect to the family of dilation, $H_{t}$ stands for the associated heat semi-group and the constant $C=C(s)$ depends on the group $\mathbb{G}$.
\end{Theoreme} 
This estimate will be a consequence of the sub-Laplacian $\mathcal{J}$ several spectral properties and we will specially use operators of type $m(\mathcal{J})$ where $m$ is a well suited Borel function (see section \ref{spectre} for the details). \\

Finally, we will prove a in very straightforward way the next theorem for a weighted functional setting: 
\begin{Theoreme}\label{PGLR}
Let $\mathbb{G}$ be a stratified Lie group and $\omega \in A_{p}$ with $1<p<+\infty$. If $f\in \dot{W}^{s_1,p}(\mathbb{G},\omega)$ and $f\in \dot{B}^{-\beta,\infty}_{\infty}(\mathbb{G})$ then
\begin{equation}\label{PGLR1}
\|f\|_{\dot{W}^{s,q}(\omega)}\leq C \|f\|_{\dot{W}^{s_1,p}(\omega)}^\theta\|f\|_{\dot{B}^{-\beta,\infty}_{\infty}}^{1-\theta}
\end{equation}
where $1<p<q<+\infty$, $\theta=p/q$, $s=\theta s_1 -(1-\theta)\beta$ and $-\beta<s<s_1$.
\end{Theoreme} 
In the proof of this theorem we will show as annouced how to pass over the Littlewood-Paley theory.\\


The plan of the article is the following: in section \ref{def} we make a short presentation of stratified Lie groups; in sections \ref{Herbata} and \ref{espa} we define maximal functions, Muckenhoupt weights and weighted functional spaces respectively; we detail the necessary results concerning spectral resolution of the sub-Laplacian in section \ref{spectre}; and finally, in section \ref{demos}, we give the proof of theorems \ref{smile00}, \ref{CHAME} and \ref{PGLR}.  

\section{Notation and preliminaries}\label{def} 

In this section we recall some basic facts about stratified Lie groups, for further information see \cite{Folland2}, \cite{Varopoulos},\cite{Stein2}  and the references given therein. \\

A \textit{homogeneous group} $\mathbb{G}$ is the data of $\mathbb{R}^{n}$ equipped with a structure of Lie group and with a family of dilations 
which are group automorphisms. We will always suppose that the origin is the identity. For \textit{dilations}, we define them by fixing integers $(a_{i})_{1\leq i\leq n}$ such that $1=a_{1}\leq... \leq a_{n}$ and by writing:  
\begin{eqnarray}\label{DefDilat}
\delta_{\alpha }:  \mathbb{R}^{n } & \longrightarrow & \mathbb{R}^{n } \label{dilat} \\
x &\longmapsto & \delta_{\alpha}[x]=(\alpha^{a_{1}}x_{1},...,\alpha^{a_{n}}x_{n})\nonumber 
\end{eqnarray} 
We will often note $\alpha x$ instead of $\delta_{\alpha}[x]$ and $\alpha$ will always indicate a strictly positive real number. \\  

The reader can be easily convinced that the Euclidean space $\mathbb{R}^{n}$ with its group structure and provided with its usual dilations 
(i.e. $a_{i}=1$, for $i=1,...,n$) is a homogeneous group. Here is another example: if $x=(x_{1}, x_{2}, x_{3})$ is an element of $\mathbb{R}^{3}$, we can fix a dilation by writing 
$\delta_{\alpha}[x]=(\alpha x_{1 }, \alpha x_{2}, \alpha^{2 } x_{3})$ for $\alpha>0$. Then, the well suited group law with respect to this dilation is given by  
\begin{equation*}
 x\cdot y=(x_{1}, x_{2}, x_{3})\cdot(y_{1}, y_{2}, y_{3})=(x_{1}+y_{1}, x_{2}+y_{2}, x_{3}+y_{3}+\frac{1}{2}(x_{1}y_{2}-y_{1}x_{2 })). 
\end{equation*} 
This triplet $(\mathbb{R}^{3}, \cdot,\delta )$ corresponds to the Heisenberg group $\mathbb{H}^{1}$ which is the first non-trivial example of a homogeneous group.\\  

The \emph{homogeneous dimension} with respect to dilation (\ref{dilat}) is given by the sum of the exponents of dilation:  
\begin{equation*} N=\sum_{1\leq i\leq n}a_{i }.  
\end{equation*} 
We observe that it is always larger than the topological dimension $n$ since integers $a_{i}$ verifies $a_{i}\geq 1$ for all $i=1,...,n$. 
For instance, in the Heisenberg group $\mathbb{H}^{1}$ we have $N=4$ and $n=3$ while in the Euclidean case these two concepts coincide. \\  

We will say that a function on $\mathbb{G}\setminus \{0\}$ is homogeneous of degree $\lambda \in \mathbb{R}$ if $f(\delta_{\alpha}[x])=\alpha^{\lambda}f(x)$ 
for all $\alpha>0$.  
In the same way, we will say that a differential operator $D$ is homogeneous of degree $\lambda$ if 
$$D(f(\delta_{\alpha}[x]))=\alpha^{\lambda}(Df)(\delta_{\alpha}[x])$$ for all $f$ in operator's domain.  
In particular, if $f$ is homogeneous of degree $\lambda$  and if $D$ is a differential operator of degree $\mu$, then $Df$ is homogeneous of degree $\lambda-\mu$. \\ 

From the point of view of measure theory, homogeneous groups behave in a traditional way since Lebesgue measure $dx$ is bi-invariant and coincides with the Haar measure. For any subset $E$ of $\mathbb{G}$ we will note its measure as $|E|$. The convolution of two functions $f$ and $g$ on $\mathbb{G}$ is defined by
\begin{equation*} 
f\ast g(x)=\int_{\mathbb{G}}f(y)g(y^{-1}\cdot x)dy=\int_{\mathbb{G}}f(x\cdot y^{-1})g(y)dy,  \quad x\in \mathbb{G}.
\end{equation*}
We also have the useful Young's inequalities:
\begin{Lemme}\label{LemYoungG} 
If $1\leq p, q, r\leq +\infty$ such that $1+\frac{1}{r}=\frac{1}{p}+\frac{1}{q}$.  
If $f\in L^{p}(\mathbb{G})$ and $g\in L^{q}(\mathbb{G})$, then $f\ast g \in L^{r}(\mathbb{G})$ and 
\begin{equation}\label{Young} 
\|f\ast g\|_{L^r}\leq \|f\|_{L^p}\|g\|_{L^q}. 
\end{equation} 
\end{Lemme} 
A proof is given in \cite{Folland2}.\\

For a homogeneous group $\mathbb{G}=(\mathbb{R}^{n}, \cdot, \delta)$ we consider now its Lie algebra $\mathfrak{g}$ whose elements can be conceived in two different ways: as \textit{left}-invariant vector fields or as \textit{right}-invariant vector fields. The left-invariant vectors fields $(X_j)_{1\leq j\leq n}$ are determined by the formula
\begin{equation}\label{loiG} 
(X_{j}f)(x)=\left.\frac{\partial f(x\cdot y)}{\partial y_{j}}\right|_{y=0}=\frac{\partial f}{\partial x_{j}}+\sum_{j<k}q^{k}_{j}(x)\frac{\partial f}{\partial x_{k}} 
\end{equation} 
where $q^{k}_{j}(x)$ is a homogeneous polynomial of degree $a_{k}-a_{j}$ and $f$ is a smooth function on $\mathbb{G}$. By this formula one deduces easily that these vectors fields are homogeneous of degree $a_{j}$:  
\begin{equation*} 
X_{j}\left(f(\alpha x)\right)=\alpha^{a_{j}}(X_{j}f)(\alpha x).  
\end{equation*} 
We will note $(Y_{j})_{1\leq j\leq n}$ the right invariant vector fields defined in a totally similar way:
$$
(Y_{j}f)(x)=\left.\frac{\partial f(y\cdot x)}{\partial y_{j}}\right|_{y=0}
$$

A homogeneous group $\mathbb{G}$ is  \emph{stratified} if its Lie algebra $\mathfrak{g}$ breaks up into a sum of linear subspaces  
$\mathfrak{g}=\bigoplus_{1\leq j\leq k } E_{j}$ such that $E_{1}$ generates the algebra $\mathfrak{g}$ and $[E_{1}, E_{j}]=E_{j+1}$ for $1\leq j < k$ and $[E_{1}, E_{k}]=\{0\}$ and $E_{k}\neq\{0\}$, but $E_{j}=\{0\}$ if $j>k$. Here $[E_{1}, E_{j}]$ indicates the subspace of $\mathfrak{g}$ generated by the elements $[U, V]=UV-VU$ with $U\in E_{1}$ and $V\in E_{j}$. The integer $k$ is called the \emph{degree} of stratification of $\mathfrak{g}$. For example, on Heisenberg group $\mathbb{H}^1$, we have $k=2$ while in the Euclidean case $k=1$.  \\

We will suppose henceforth that $\mathbb{G}$ is \textbf{stratified}. Within this framework, if we fix the vectors fields $X_{1},...,X_{m}$ such that $a_{1}=a_{2}=\ldots=a_{m}=1$ $(m<n)$, then the family $(X_{j})_{1\leq j\leq m}$ is a base of $E_{1}$ and generates the Lie algebra of $\mathfrak{g}$, which is precisely the H\"ormander's condition (see \cite{Folland2} and \cite{Varopoulos}).\\

To the family $(X_{j})_{1\leq j\leq m}$ is associated the Carnot-Carath\'eodory distance $d$ which is left-invariant and compatible with the topology on $\mathbb{G}$ (see \cite{Varopoulos} for more details). For any $x\in \mathbb{G}$ we will note $|x|=d(x,e)$ and for $r>0$ we form the balls by writing $B(x,r)=\{y\in \mathbb{G}: d(x,y)<r\}$.\\

The main tools of this paper depends on the properties of the gradient, the sub-Laplacian and the associated heat kernel. Before introducing them, we make here three remarks on general vectors fields $X_{j}$ and $Y_{j}$. Let us fix some notation. For any multi-index $I=(i_{1},...,i_{n})\in \mathbb{N}^{n}$, one defines $X^{I}$ by $X^{I}=X_{1}^{i_{1}}\dots X_{n}^{i_{n}}$ and $Y^{I}$ by $Y^{I}=Y_{1}^{i_{1}}\dots Y_{n}^{i_{n}}$. We note $|I|= i_{1}+\ldots+i_{n}$ the order of the derivation $X^I$ or $Y^I$ and $d(I)=a_{1}i_{1}+\ldots+a_{n}i_{n}$ the homogeneous degree of this one. \\  

Firstly, for $\varphi, \psi \in \mathcal{C}^{\infty}_{0}(\mathbb{G})$ we have the equality 
\begin{equation*} 
\int_{\mathbb{G}}\varphi(x)(X^{I}\psi)(x)dx=(-1)^{|I|}\int_{\mathbb{G}}(X^{I}\varphi)(x)\psi(x)dx.  
\end{equation*} 
Secondly, interaction of operators $X^{I}$ and $Y^{I}$ with convolutions is clarified by the following identities:  
\begin{equation}\label{izder}
X^{I}(f*g)=f*(X^{I}g), \qquad Y^{I}(f*g)=(Y^{I}f)*g, \qquad (X^{I}f)*g=f*(Y^{I}g). 
\end{equation} 
Finally, one will say that a function $f \in \mathcal{C}^{\infty}(\mathbb{G})$ belongs to the Schwartz class $\mathcal{S}(\mathbb{G})$ if the following semi-norms are bounded for all $k\in \mathbb{N}$ and any multi-index $I$:  $N_{k,I}(f)=\underset{x\in \mathbb{G}}{\sup } \, (1+|x|)^{k }|X^{I}f(x)|$. 
\begin{Remarque}
\emph{To characterize the Schwartz class $\mathcal{S}(\mathbb{G})$ we can replace vector fields $X^I$ in the semi-norms $N_{k,I}$ above by right-invariant vector fields $Y^I$.}
\end{Remarque}
For a proof of these facts and for further details see \cite{Folland2} and \cite{Furioli2}.\\

We define now the \textit{gradient} on $\mathbb{G}$ from vectors fields of homogeneity degree equal to one by fixing $\nabla = (X_{1},...,X_{m })$.  This operator is of course left invariant and homogeneous of degree $1$. The length of the gradient is given by the formula $|\nabla f|= \left((X_{1}f)^{2}+... +(X_{m}f)^{2 } \right)^{1/2}$.\\
 
Let us notice that there is not a single way to build a \textit{sub-Laplacian}, see for example \cite{Furioli2} and \cite{Chame}. In this article, we will work with the following sub-Laplacian:
\begin{equation}\label{laplaciano} 
\mathcal{J}=\nabla^{*}\nabla=-\sum_{j=1}^{m}X^{2}_{j} 
\end{equation} 
which is a positive self-adjoint, hypo-elliptic operator (since $(X_j)_{1\leq j\leq m}$ satisfies the H\"ormander's condition), having as domain of definition $L^2(\mathbb{G})$. Its associated \textit{heat operator} on $\mathbb{G}\times]0, +\infty[$ is given by $\partial_{t}+\mathcal{J}$.\\

We recall now some well-known properties of this operator. 
 
\begin{Theoreme}\label{fol} There exists a unique family of continuous linear operators $(H_{t})_{t>0}$ defined on $L^{1}+L^{\infty}(\mathbb{G})$ with the semi-group property $H_{t+s}=H_{t}H_{s}$ for all $t, s>0$ and $H_{0}=Id$, such that: \\  
\begin{enumerate} 
\item[1)] the sub-Laplacian $\mathcal{J}$ is the infinitesimal generator of the semi-group  $H_{t}=e^{-t\mathcal{J}}$; \\  
\item[2)] $H_{t}$ is a contraction operator on $L^{p}(\mathbb{G})$ for $1\leq p\leq +\infty$ and for $t>0$; \\  
\item[3)] the semi-group $H_t$ admits a convolution kernel $H_{t}f=f\ast h_{t}$ where $h_{t}(x)=h(x, t) \in \mathcal{C}^{\infty}(\mathbb{G}\times]0, +\infty[)$ is the heat kernel which satisfies the following points:  
\begin{enumerate} 
\item $(\partial_{t}+\mathcal{J})h_{t}=0$ on $\mathbb{G}\times]0, +\infty[$, \\ 
\item $h(x, t)=h(x^{-1}, t)$, $h(x, t)\geq 0$ and $\int_{\mathbb{G}}h(x, t)dx=1$, \\ 
\item $h_{t}$ has the semi-group property:  $h_{t}\ast h_{s}=h_{t+s}$ for $t, s>0$, \\ 
\item $h(\delta_{\alpha}[x], \alpha^{2}t)=\alpha^{-N}h(x, t)$, \\ 
\item For every $t>0$, $x\mapsto h(x, t)$ belong to the Schwartz class in $\mathbb{G}$. \\  
\end{enumerate} 
\item[4)] $ \|H_{t}f-f \|_ {L^p}\to 0$ if $t\to 0$ for $f \in L^{p}(\mathbb{G})$ and $1\leq p < +\infty$; \\  
\item[5)] If $f\in L^{p}(\mathbb{G})$, $1\leq p\leq +\infty$, then the function $u(x, t)=H_{t}f(x)\in \mathcal{C}^{\infty}(\mathbb{G}\times \mathbb{R}^{+})$ is a solution of the heat equation: 
$$ \qquad\left \lbrace\begin{array}{l}(\frac{\partial}{\partial t}+\mathcal{J})u(x, t)=0\quad \mbox{for } \quad x\in \mathbb{G } \quad \mbox{and}\quad t>0\, ;\\[8mm] 
u(x, 0)=f(x)\qquad\quad \mbox{for } \quad x\in \mathbb{G }.\\  
\end{array}\right. $$ 
\end{enumerate} 
\end{Theoreme} 

For a detailed proof of these and other important facts concerning the heat semi-group see \cite{Folland2} and \cite{Saka}.\\

To close this section we recall the definition of the sub-Laplacian fractional powers $\mathcal{J}^s$ with $s>0$.\\ We write:
$$\mathcal{J}^s f(x)=\underset{\varepsilon \to 0}{\lim}\frac{1}{\Gamma(k-s)}\int_{\varepsilon}^{+\infty}t^{k-s-1}\mathcal{J}^k H_tf(x)dt$$
for all $f\in \mathcal{C}^{\infty}(\mathbb{G})$ with $k$ the smallest integer greater than $s$.

\section{Maximal functions and Muckenhoupt weights}\label{Herbata}
There are several ways of defining maximal functions in stratified Lie groups and our principal reference is \cite{Folland2}. In this article we will mainly work with the following function  
\begin{Definition} Let $f\in \mathcal{S}'(\mathbb{G})$ and $\varphi\in\mathcal{S}(\mathbb{G})$.  
The maximal function $\mathcal{M}_{\varphi}$ is given by the expression 
\begin{equation*} 
\mathcal{M}_{\varphi}f(x)=\underset{0<t<+\infty}{\sup}\{|f\ast\varphi_{t}(x)|\} 
\end{equation*} with $\varphi_{t}(x)=t^{-N/2}\varphi(t^{-1/2}x)$.  
\end{Definition} 
This definition still has a sense if $f$ and $\varphi$ are two distributions such that $(x, t)\longmapsto f\ast\varphi_{t}(x)$ is a continuous function on $\mathbb{G}\times]0, +\infty[$:  for example, if $f\in L^{p}(\mathbb{G})$ and $\varphi\in L^{q}(\mathbb{G})$ where $1\leq p\leq +\infty$ and $1/p+1/q=1$. \\ 

An important special case is given by the Hardy-Littlewood function which consists in taking as function $\varphi$ the characteristic function of the unit ball:  
$$\mathcal{M}_{B}f(x)=\underset{B \ni x}{\sup } \;\frac{1}{|B|}\int_{B }|f(y)|dy.$$ 
The next lemma explain the relationship between these maximal functions.
\begin{Lemme}\label{Led11} 
Let $\varphi$ a function on $\mathbb{G}$ such that $|\varphi(x)|\leq C(1+|x|)^{-N-\varepsilon}$ for some $\varepsilon>0$, then
\begin{equation}\label{Wilma} 
\mathcal{M}_{\varphi}f(x)\leq C \mathcal{M}_{B}f(x). 
\end{equation} 
\end{Lemme} 
We will use this property in the sequel and we request the reader to consult the proof in \cite{Folland2}. \\  

The reader can consult \cite{Folland2}, \cite{Grafakos} and \cite{Garcia} for a more detailed study of these important functions. 
For our part, we will be interested in the relationship existing between these functions and weights. A \emph{weight} $\omega$ is, in a very general way, a locally integrable function on $\mathbb{G}$ with values in $]0,+\infty[$.  
For a given weight $\omega$ and a measurable set $E\subset\mathbb{G}$ we use the following notation:  
\begin{equation*} 
\omega(E)=\int_{E}\omega(x)dx.  
\end{equation*} 
We will define thus, for $1\leq p<+\infty$, weighted Lebesgue spaces by the norm 
\begin{equation}\label{pico1}
 \|f\|_ {L^{p}(\omega)}=\left(\int_{\mathbb{G}}|f(x)|^{p}\omega(x)dx\right)^{1/p }
\end{equation}
Historically, the characterization of Muckenhoupt weights comes from the following problem: for a fixed $p\in ]1, +\infty[$ we want to know for which functions 
$\omega$ one has the strong estimate
\begin{equation}\label{MaximalProper}
\int_{\mathbb{G}}\mathcal{M}_{B}f(x)^{p}\omega(x)dx\leq C \int_{\mathbb{G}}|f(x)|^{p}\omega(x)dx \qquad (f\in L^{p}(\mathbb{G},\omega)).  
\end{equation} 
It follows the condition below and the next definition (see \cite{Grafakos}): 
\begin{equation}\label{brown} 
\underset{B}{\sup}\left(\frac{1}{|B|}\int_{B}\omega(x)dx\right)\left(\frac{1}{|B|}\int_{B}\omega(x)^{-\frac{1}{p-1}}dx\right)^{p-1}<+\infty. 
\end{equation} 

\begin{Definition} 
Let $\mathbb{G}$ a stratified Lie group and let $1<p<+\infty$.  We will say that a weight $\omega$ belongs to the Muckenhoupt class $A_{p}$ if it satisfies condition (\ref{brown}). Moreover, we will define weights in the class $A_{1}$ by:  
\begin{equation}\label{soledad} 
\mathcal{M}_{B}\;\omega(x)\leq C \;\omega(x) \qquad (\forall x\in \mathbb{G}).  
\end{equation} 
\end{Definition}

Here are some traditional examples:  the trivial weight $\omega(x) \equiv 1$ for all $x\in \mathbb{G}$ is a $A_p$ weight for $1\leq p<+\infty$ and the function $|x|^{\alpha}$ is in $A_{p}$ if and only if  $-N<\alpha<N(p-1)$, where $N$ is the homogeneous dimension. For $p=1$, the function $|x|^{\alpha}$ belongs to $A_{1}$ if and only if $-N<\alpha\leq 0$. \\  

Let us finally say that we have following inclusion:
\begin{Proposition}\label{soledad1bis} 
If $1< p<q<+\infty$, then $A_{1}\subset A_{p}\subset A_{q}$.
\end{Proposition} 
We request the reader to consult the proof of this result in \cite{Folland2} or \cite{Grafakos}.  
\section{Weigthed spaces}\label{espa} 

We give in this section the precise definition of weighted functional spaces involved in theorems \ref{smile00}, \ref{CHAME} and \ref{PGLR}. In a general way, given a norm $\|\cdot\|_{X(\omega)}$, we will define the corresponding weighted functional space $X(\mathbb{G}, \omega)$ by $\{f\in \mathcal{S}'(\mathbb{G}): \|f\|_{X(\omega)}<+\infty\}$ where $\omega$ is a Muckenhoupt weight belonging to a certain class $A_p$.\\

\begin{enumerate}
\item[$\bullet$] \textbf{Lebesgue spaces} $L^p(\mathbb{G}, \omega)$. We have already considered how to define weigthed Lebesgue spaces with the formula (\ref{pico1}). Let us notice that we also have a characterization with the distribution function:
\begin{equation*} 
\|f\|^{p}_{L^{p}(\omega)}=\int_{0}^{+\infty}p\sigma^{p-1 } \omega(\{x\in \mathbb{G}:|f(x)|> \sigma\})d\sigma.  
\end{equation*} 
\item[$\bullet$] \textbf{weak-$L^{p}$ spaces} or Lorentz spaces $L^{p,\infty}(\mathbb{G}, \omega)$. We define them by
\begin{equation*}
\|f\|_ { L^{p, \infty}(\omega)}=\underset{\sigma>0}{\sup}\{\sigma \; \omega(\{x\in \mathbb{G}:|f(x)|> \sigma\})^{1/p}\}.
\end{equation*} 
\item[$\bullet$]\textbf{Sobolev spaces} $\dot{W}^{s,p}(\mathbb{G}, \omega)$. For $\omega\in A_{p}$ we write:  
\begin{equation*} 
\|f\|_ {\dot{W}^{s, p}(\omega)} =\|\mathcal{J}^{s/2}f\|_{L^{p}(\omega)} \qquad (1<p<+\infty) 
\end{equation*}
and when $p=s=1$ we will note
\begin{equation*}
\|f\|_ {\dot{W}^{1,1}(\omega)} =\|\nabla f\|_{L^{1}(\omega)}.
\end{equation*} 
\item[$\bullet$]\textbf{weak Sobolev spaces} $\dot{W}^{s,p}_\infty(\mathbb{G}, \omega)$: 
\begin{equation} \label{DefDebilSoboLe}
\|f\|_ {\dot{W}^{s, p}(\omega)} =\|\mathcal{J}^{s/2}f\|_{L^{p,\infty}(\omega)} \qquad (1<p<+\infty) 
\end{equation}
\item[$\bullet$] \textbf{Besov spaces} $\dot{B}^{s,q}_p(\mathbb{G}, \omega)$. We define them in the following way: 
\begin{equation*} 
\|f\|_{\dot{B}^{s, q}_{p}(\omega)}=\left[\int_{0}^{+\infty}t^{(m-s/2)q}\left \|\frac{\partial^{m}H_{t}f}{\partial t^{m}}(\cdot)\right\|^{q}_{L^{p}(\omega)}\frac{dt}{t}
\right]^{1/q} 
\end{equation*} 
for $1\leq p, q\leq +\infty, s>0$ and $m$ an integer such that $m>s/2$. \\ 
 
Finally, for Besov spaces of indices $(-\beta, \infty, \infty)$ which appear in all the Improved Sobolev inequalities we have:
\begin{equation}\label{BesovChaleurGhomo}
\|f\|_{\dot{B}^{-\beta,\infty}_{\infty}}=\underset{t>0}{\sup}\;\;t^{\beta/2 } \|H_{t}f\|_ {L^\infty} 
\end{equation} 
\end{enumerate}
\section{Spectral resolution of the sub-Laplacian}\label{spectre}

The use in this article of spectral resolution for the sub-Laplacian consists roughly in expressing this operator by the formula $\mathcal{J}=\int_{0}^{+\infty}\lambda \;dE_{\lambda}$ and, by means of this characterization, build a family of new operators $m(\mathcal{J})$ associated to a Borel function $m$. This kind of operators have some nice properties as shown in the next propositions.
\begin{Proposition} 
If $m$ is a bounded Borel function on $]0,+\infty[$ then the operator $m(\mathcal{J})$ fixed by
\begin{equation}\label{little} 
m(\mathcal{J})=\int_{0}^{+\infty}m(\lambda) \;dE_{\lambda }, 
\end{equation} 
is bounded on $L^{2}(\mathbb{G})$  and admits a convolution kernel $M$ i.e.: $m(\mathcal{J})(f)=f\ast M \qquad (\forall f\in L^{2}(\mathbb{G}))$. 
\end{Proposition} 
See \cite{Folland2} and the references given therein for a proof. For our purposes, it will be particularly interesting to combine this result with the structure of dilation:  
\begin{Lemme}\label{madera} 
Let $m$ be a bounded function on $]0,+\infty[$ and let $M$ be the kernel of the operator $m(\mathcal{J})$.  
Then, for all $t>0$ we can build a bounded operator on $L^{2}(\mathbb{G})$ by writing $m_{t}(\mathcal{J})=m(t\mathcal{J})$ with an associated kernel given by 
\begin{equation*} 
M_{t}(x)=t^{-N/2}M(t^{-1/2}x).  
\end{equation*} 
\end{Lemme} 

Following \cite{Hulanicki} and \cite{Furioli2} we can improve the conclusion of the above proposition. Let $k\in\mathbb{N}$ and $m$ be a function of class $\mathcal{C}^{k}(\mathbb{R}^{+})$, we write 
\begin{equation*} 
\|m \|_ {(k)}=\underset{\underset{\lambda>0}{1\leq r\leq k}}{\sup } (1+\lambda)^{k }|m^{(r)}(\lambda)|.
\end{equation*}
This formula gives us a necessary condition to obtain certain properties of the operators defined by (\ref{little}):
\begin{Proposition}\label{toge}
Let $\alpha\in \mathbb{N}$, $I=(i_{1},...,i_{n})$ be a multi-index and $p\in [1, +\infty]$.  
There is a constant $C>0$ and an integer $k$ such that, for any function $m\in \mathcal{C}^{k}(\mathbb{R}^{+})$ with $ \|m\|_{(k)}<+\infty$, the kernel $M_{t}$ associated to the operator $m(t\mathcal{J})$, $t>0$, satisfies 
 \begin{equation*} 
\|(1+|\cdot|)^{\alpha}X^{I}M_{t}(\cdot)\|_{L^p}\leq C(1+\sqrt{t})^{\alpha}t^{-(\frac{N}{2p'}+\frac{d(I)}{2})}\|m\|_{(k)}. 
\end{equation*} 
where $\frac{1}{p}+\frac{1}{p'}=1$.
\end{Proposition} 
\begin{Corollaire}\label{black} Let $t>0$.
\begin{enumerate}
\item[1)] Let $m$ be the restriction on $\mathbb{R}^{+}$ of a function defined on $\mathcal{S}(\mathbb{R})$.  Then, the kernel $M$ 
of the operator $m(\mathcal{J})$ is in $\mathcal{S}(\mathbb{G})$.  
\item[2)] If $m$ is as above; and if it is vanishing at all orders near of the origin, then the kernel $M$ belongs to the 
space $\mathcal{S}_{0}(\mathbb{G})$ formed by the functions of the Schwartz class which every moment is null.  
\end{enumerate} 
\end{Corollaire}
For more details and proofs see \cite{Folland2}, \cite{Furioli2} and \cite{Hulanicki}.
\section{Improved Sobolev Inequalities on stratified groups: the proofs}\label{demos} 

As said in the introduction, inequalities given in theorem \ref{smile00} depends on the theorem \ref{CHAME}. We will thus begin proving this result in the following lines and we will continue our study by treating separately weak inequalities (\ref{smile1}) and strong inequalities (\ref{smile0}).
\subsection{The modified pseudo-inequality of Poincar\'e}\label{ISPGH} 

Under the hypothesis of theorem \ref{CHAME}, we have to prove the inequality  
\begin{equation*} 
\|\mathcal{J}^{s/2}f-H_{t}\mathcal{J}^{s/2}f\|_{L^{1}(\omega)}\leq C \; t^{\frac{1-s}{2}}\|\nabla f\|_{ L^{1}(\omega)}.  
\end{equation*} 
To begin the proof, we observe that the following identity occurs: 
\begin{equation*} 
(\mathcal{J}^{s/2}f-H_{t}\mathcal{J}^{s/2}f)(x)=\left(\int_{0}^{+\infty}m(t\lambda)dE_{\lambda}\right)t^{1-s/2}\mathcal{J} f(x), 
\end{equation*} 
where we noted $m(\lambda)=\lambda^{s/2-1}(1-e^{-\lambda }$) for $\lambda>0$, note that $m$ is a bounded function which tends to $0$ at infinity since $s/2-1<0$. We break up this function by writing:  
$$m(\lambda)=m_{0}(\lambda)+m_{1}(\lambda)=m(\lambda)\theta_{0}(\lambda)+m(\lambda)\theta_{1}(\lambda)$$ 
where we chose the auxiliary functions $\theta_{0}(\lambda), \theta_{1}(\lambda)\in \mathcal{C}^{\infty}(\mathbb{R}^{+})$ defined by:  
\begin{eqnarray*} 
\bullet\quad \theta_{0}(\lambda) = 1 \quad\mbox{on } \quad]0, 1/2 ] \quad\mbox{and}\quad 0 \quad \mbox{on}\quad]1, +\infty[, \\[5mm ] 
\bullet\quad \theta_{1}(\lambda) = 0 \quad\mbox{on} \quad]0, 1/2 ] \quad\mbox{and}\quad 1 \quad \mbox{on }\quad]1, +\infty[, 
\end{eqnarray*} 
so that $\theta_{0}(\lambda)+\theta_{1}(\lambda)\equiv 1$. Then, we obtain the formula: 
\begin{equation*} 
(\mathcal{J}^{s/2}f-H_{t}\mathcal{J}^{s/2}f)(x)=\left(\int_{0}^{+\infty}m_{0}(t\lambda)dE_{\lambda}\right)t^{1-s/2}\mathcal{J}f(x)+
\left(\int_{0}^{+\infty}m_{1}(t\lambda)dE_{\lambda}\right)t^{1-s/2}\mathcal{J}f(x).  
\end{equation*} 
If we note $M^{(i)}_{t}$ the kernel of the operator fixed by $\int_{0}^{+\infty}m_{i}(t\lambda)dE_{\lambda}$ for $i=0,1$, we have:  
\begin{equation*} 
(\mathcal{J}^{s/2}f-H_{t}\mathcal{J}^{s/2}f)(x)=t^{1-s/2}\mathcal{J}f\ast M^{(0)}_{t}(x)+t^{1-s/2}\mathcal{J}f\ast M^{(1)}_{t}(x).  
\end{equation*} 
We now multiply the above equality by a weight $\omega\in A_1$ to obtain the inequality
\begin{equation}\label{poids1} 
\int_{\mathbb{G}}\left|\mathcal{J}^{s/2}f-H_{t}\mathcal{J}^{s/2}\right|\omega(x)dx \leq \int_{\mathbb{G}}\left|t^{1-s/2}\mathcal{J}f\ast M^{(0)}_{t}(x)\right|\omega(x)dx +\int_{\mathbb{G}}\left|t^{1-s/2}\mathcal{J}f\ast M^{(1)}_{t}(x)\right|\omega(x)dx.
\end{equation} 
We will now estimate the right side of the above inequality by the two following propositions:  
\begin{Proposition}\label{Sir} 
For the first integral in the right-hand side of (\ref{poids1}) we have the inequality:  
\begin{equation*} 
\int_{\mathbb{G}}\left|t^{1-s/2}\mathcal{J}f\ast M^{(0)}_{t}(x)\right|\omega(x)dx\leq Ct^{\frac{1-s}{2}}\|\nabla f\|_{L^{1}(\omega)} 
\end{equation*}
\end{Proposition} 
\emph{\textbf{Proof}}. The function $m_{0}$ is the restriction on $\mathbb{R}^{+}$ of a function belonging to the Schwartz class. This function satisfies the assumptions of corollary \ref{black} which we apply after having noticed the identity
\begin{equation*} 
I=\int_{\mathbb{G}}\left|t^{1-s/2}\mathcal{J}f\ast M^{(0)}_{t}(x)\right|\omega(x)dx=
\int_{\mathbb{G}}\left|t^{1-s/2}\nabla f\ast\tilde{\nabla}M^{(0)}_{t}(x)\right|\omega(x)dx 
\end{equation*} 
where we noted $\tilde{\nabla}$ the gradient formed by vectors fields $(Y_{j})_{1\leq j\leq m}$. We have then
\begin{equation*} 
I\leq \int_{\mathbb{G}}\int_{\mathbb{G}}t^{1-s/2}|\nabla f(y)||\tilde{\nabla} M^{(0)}_{t}(y^{-1}\cdot x)|\omega(x)dxdy.  
\end{equation*} 
By corollary \ref{black}, one has $M^{(0)}_{t}\in \mathcal{S}(\mathbb{G})$ and, since $M_{t}^{(0)}(x)=t^{-N/2}M^{(0)}(t^{-1/2}x)$, we can write 
$$K_{t}(x)=t^{1/2 }|\tilde{\nabla}M^{(0)}_{t}(x^{-1 })|\in L^{1}(\mathbb{G}).$$ 
One obtains
\begin{equation*} 
I\leq \int_{\mathbb{G}}\int_{\mathbb{G}}t^{\frac{1-s}{2}}|\nabla f(y)|\, K_{t}(x\cdot y^{-1})\omega(x)dxdy=\int_{\mathbb{G}}t^{\frac{1-s}{2}}|\nabla f(y)|\; 
\omega\ast K_{t}(y)dy.  
\end{equation*} 
By definition of maximal functions and by the estimate (\ref{Wilma}), we have the inequality:  
$$\underset{t>0}{\sup } \;\omega\ast K_{t}(y)\leq C\, \left(\mathcal{M}_{B}\;\omega\right)(y), $$ 
hence,   
\begin{equation*}
 I\leq Ct^{\frac{1-s}{2}}\int_{\mathbb{G}}|\nabla f(y)|\mathcal{M}_{B}\, \omega(y)dy.  
\end{equation*} 
It remains to notice that, by assumption, $\omega \in A_{1}$ if and only if $(\mathcal{M}_{B}\, \omega)(\cdot)\leq C \;\omega(\cdot).$ 
We obtain then the desired estimation:  
\begin{equation*} 
\int_{\mathbb{G}}\left|t^{1-s/2}\mathcal{J}f\ast M^{(0)}_{t}(x)\right|\omega(x)dx\leq Ct^{\frac{1-s}{2}}\int_{\mathbb{G}}|\nabla f(y)|\omega(y)dy. 
\end{equation*} 
\begin{flushright}{$\blacksquare$}\end{flushright}
\begin{Proposition}\label{PropoChame2}
For the last integral of (\ref{poids1}) we have the inequality
\begin{equation*} 
\int_{\mathbb{G}}\left|t^{1-s/2}\mathcal{J}f\ast M^{(1)}_{t}(x)\right|\omega(x)dx\leq Ct^{\frac{1-s}{2}}\|\nabla f\|_{L^{1}(\omega)} 
\end{equation*}
\end{Proposition} 
\emph{\textbf{Proof}}. Here, it is necessary to make an additional step. We cut out the function $m_{1}$ in the following way:  
\begin{equation*}
m_{1}(\lambda)=\left(\frac{1-e^{-\lambda}}{\lambda}\right)\theta_{1}(\lambda)=m_{a}(\lambda)-m_{b}(\lambda) 
\end{equation*} 
where $m_{a}(\lambda)=\frac{1}{\lambda}\theta_{1}(\lambda)$ and $m_{b}(\lambda)=\frac{e^{-\lambda}}{\lambda}\theta_{1}(\lambda)$. We will note $M^{(a)}_{t}$ and $M^{(b)}_{t}$ the associated kernels of these two operators. We obtain thus the estimate 
\begin{equation}\label{Rolling1} 
\int_{\mathbb{G}}\left|t^{1-s/2}\mathcal{J}f\ast M^{(1)}_{t}(x)\right|\omega(x)dx\leq \int_{\mathbb{G}}\left|t^{1-s/2}\mathcal{J}f\ast M^{(a)}_{t}(x)\right|\omega(x)dx+ 
\int_{\mathbb{G}}\left|t^{1-s/2}\mathcal{J}f\ast M^{(b)}_{t}(x)\right|\omega(x)dx 
\end{equation} 
Observe that $m_{b}\in \mathcal{S}(\mathbb{R}^{+})$ and then $M^{(b)}_{t}\in \mathcal{S}(\mathbb{G})$.  We have the next lemma for the last integral in (\ref{Rolling1}).
\begin{Lemme}\label{BishopLem0}
\begin{equation*} 
\int_{\mathbb{G}}\left|t^{1-s/2}\mathcal{J}f\ast M^{(b)}_{t}(x)\right|\omega(x)dx\leq Ct^{\frac{1-s}{2}}\|\nabla f\|_{L^{1}(\omega)}.
\end{equation*} 
\end{Lemme} 
\emph{\textbf{Proof}}. The proof is straightforward and follows the same steps as those of the preceding proposition \ref{Sir}.  
\begin{flushright}{$\blacksquare$}\end{flushright} 
We treat the other part of (\ref{Rolling1}) with the following lemma:  
\begin{Lemme} \label{BishopLem}
\begin{equation}\label{Bishop} 
\int_{\mathbb{G}}\left|t^{1-s/2}\mathcal{J}f\ast M^{(a)}_{t}(x)\right|\omega(x)dx\leq Ct^{\frac{1-s}{2}}\|\nabla f\|_{L^{1}(\omega)} 
\end{equation} \end{Lemme} 
\emph{\textbf{Proof}}. We consider the auxiliary function 
$$\psi(\lambda)=\theta_{0}(\lambda/2)-\theta_{0}(\lambda)=\theta_{1}(\lambda)-\theta_{1}(\lambda/2)$$ 
in order to obtain the identity $$\sum_{j=0}^{+\infty}\psi(2^{-j}\lambda)=\theta_{1}(\lambda).$$  
We have then 
$$m_{a}(t\lambda)=\frac{1}{t\lambda}\sum_{j=0}^{+\infty}\psi(2^{-j}t\lambda)=\sum_{j=0}^{+\infty}2^{-j}\tilde{\psi}(2^{-j}t\lambda)$$ 
where $\tilde{\psi}(\lambda)=\frac{\psi(\lambda)}{\lambda}$ is a function in $\mathcal{C}^{\infty}_{0}(\mathbb{R}^{+})$. By corollary \ref{black}, the kernel $\tilde{K}$ associated with the function $\tilde{\psi}$ belongs to $\mathcal{S}_{0}(\mathbb{G})$. Then, from the point of view of operators, one has:  
\begin{equation}\label{Stone} 
M^{(a)}_{t}(x)=\sum^{+\infty}_{j=0}2^{-j}\tilde{K}_{j, t}(x) 
\end{equation} 
where $\tilde{K}_{j,t}(x)=2^{N/2}t^{-N/2}\tilde{K}(2^{j/2}t^{-1/2}x)$. With formula (\ref{Stone}) we return to the left side of (\ref{Bishop}):  
$$\int_{\mathbb{G}}\left|t^{1-s/2}\mathcal{J}f\ast M^{(a)}_{t}(x)\right|\omega(x)dx\leq\sum^{+\infty}_{j=0}2^{-j}\int_{\mathbb{G}}\left|t^{1-s/2}\mathcal{J}f\ast 
\tilde{K}_{j,t}(x)\right|\omega(x)dx.$$
Using the sub-Laplacian definition and vector fields properties, we have
\begin{equation}\label{Chicago}\int_{\mathbb{G}}\left|t^{1-s/2}\mathcal{J}f\ast M^{(a)}_{t}(x)\right|\omega(x)dx\leq \sum^{+\infty}_{j=0}2^{-j}t^{1-s/2}\int_{\mathbb{G}}
\int_{\mathbb{G } }|\nabla f(y)| |\tilde{\nabla}\tilde{K}_{j, t}(y^{-1}\cdot x)|\omega(x)dxdy.
\end{equation} 
We note this time $K_{j, t}(x)=2^{-j/2}t^{1/2}|\tilde{\nabla}\tilde{K}_{j, k}(x^{-1})|$ to obtain the following formula for the right side of (\ref{Chicago}):
$$\sum^{+\infty}_{j=0}2^{-j/2}t^{\frac{1-s}{2}}\int_{\mathbb{G}}\int_{\mathbb{G}}|\nabla f(y)|K_{j,t}(x\cdot y^{-1})\omega(x)dxdy = 
\sum^{+\infty}_{j=0}2^{-j/2}t^{\frac{1-s}{2}}\int_{\mathbb{G}}|\nabla f(y)|\;\omega\ast K_{j, t}(y)dy.$$ 
It remains to apply the same arguments used in proposition \ref{Sir}, namely the assumption $\omega \in A_{1}$ and, for $K_{j,t}$, the estimations $\underset{j, t>0}{\sup } \;\omega \ast K_{j, t}(y)\leq C\, (\mathcal{M}_{B } \; \omega)(y)\leq C\, \omega(y).$ 
Then, we finally get the inequality  
\begin{equation*} 
\int_{\mathbb{G}}\left|t^{1-s/2}\mathcal{J}f\ast M^{(a)}_{t}(x)\right|\omega(x)dx\leq C\, t^{\frac{1-s}{2}}\sum_{j=0}^{+\infty}2^{-j/2}\int_{\mathbb{G}}|\nabla f|
(y)\omega(y)dy = C\, t^{\frac{1-s}{2} } \|\nabla f\|_{L^{1}(\omega)}.  
\end{equation*} 
Which ends the proof of the lemma \ref{BishopLem}.
\begin{flushright}{$\blacksquare$}\end{flushright} 
With these two last lemmas we conclude the proof of the proposition \ref{PropoChame2}. Now, getting back to the formula (\ref{poids1}), with propositions \ref{Sir} and \ref{PropoChame2} we finally finish the proof of theorem \ref{CHAME}.
\begin{flushright}{$\blacksquare$}\end{flushright} 
\subsection{Weak inequalities}

To begin the proof notice that operator $\mathcal{J}^{s/2}$ carries out an isomorphism between the spaces $\dot{B}^{-\beta, \infty}_{\infty}(\mathbb{G})$ and $\dot{B}^{-\beta-s, \infty}_{\infty}(\mathbb{G})$ (see \cite{Saka}). Thus inequality (\ref{smile1}) rewrites as:  
\begin{equation}\label{ForIsoLap}
\|\mathcal{J}^{s/2}f\|_{L^{q, \infty}(\omega)}\leq C \|\nabla f\|_{L^{1}(\omega)}^{\theta}\|\mathcal{J}^{s/2}f\|_{\dot{B}^{-\beta-s, \infty}_{\infty}}^{1-\theta} 
\end{equation} 
By homogeneity, we can suppose that the norm $\|\mathcal{J}^{s/2}f \|_{\dot{B}^{-\beta-s,\infty}_{\infty}}$ is bounded by $1$; then we have to show 
\begin{equation}\label{CHAMEFAIBLE2} 
\|\mathcal{J}^{s/2}f \|_{L^{q, \infty}(\omega)}\leq C \|\nabla f \|_{L^{1}(\omega)}^{\theta}. 
\end{equation} 
We have thus to evaluate the expression $\omega\left(\{x\in \mathbb{G}:|\mathcal{J}^{s/2}f(x)|> 2\alpha\}\right)$ for all $\alpha>0$. If we use the thermic definition of the Besov space (\ref{BesovChaleurGhomo}), we have 
$$\|\mathcal{J}^{s/2}f\|_{\dot{B}^{-\beta-s,\infty}_{\infty}}\leq 1 \iff \underset{t>0}{\sup}\left\{t^{\frac{\beta+s}{2}}\|H_{t}\mathcal{J}^{s/2}f\|_{L^\infty}\right\} \leq 1.$$ 
But, if one fixes $t_{\alpha}=\alpha^{-\left(\frac{2}{\beta+s}\right)}$, we obtain $\|H_{t_{\alpha}}\mathcal{J}^{s/2}f \|_{L^\infty}\leq \alpha $. 
Note also that with the definition of parameter $\beta$ one has $t_{\alpha}=\alpha^{-\frac{2(q-1)}{(1-s)}}$. Therefore, since we have the following set inclusion
$$\left\{x\in \mathbb{G}: |\mathcal{J}^{s/2}f(x)|> 2\alpha\right\}\subset \left\{x\in \mathbb{G}: |\mathcal{J}^{s/2}f(x)-H_{t_{\alpha}}\mathcal{J}^{s/2}f(x)|> \alpha\right\}, $$ 
the Tchebytchev inequality implies 
\begin{equation*} 
\alpha^{q}\omega\left(\{x\in \mathbb{G}: |\mathcal{J}^{s/2}f(x)|> 2\alpha\}\right)\leq \alpha^{q-1}\int_{\mathbb{G } }|\mathcal{J}^{s/2}f(x)-H_{t_{\alpha}}\mathcal{J}^{s/2}f(x)|\omega(x)dx.  
\end{equation*} 
At this point, we use the modified Poincar\'e pseudo-inequality, given by theorem \ref{CHAME}, to estimate the right side of the preceding inequality:  
\begin{equation}\label{T2g} 
\alpha^{q}\omega\left(\{x\in \mathbb{G}:|\mathcal{J}^{s/2}f(x)|> 2\alpha\}\right)\leq C \alpha^{q-1}\;t_{\alpha}^{\frac{1-s}{2}}\int_{\mathbb{G}}|\nabla f(x)|\omega(x)dx.  
\end{equation} 
But, by the choice of $t_{\alpha}$, one has $\alpha^{q-1}\alpha^{-\frac{2(q-1)}{(1-s)}\frac{(1-s)}{2}}=1 $. Then (\ref{T2g}) implies the inequality 
$$\qquad\alpha^{q}\omega\left(\{x\in \mathbb{G}:|\mathcal{J}^{s/2}f(x)|> 2\alpha\}\right)\leq C \|\nabla f \|_ { L^{1}(\omega) } \;; $$  
and, finally, using definition (\ref{DefDebilSoboLe}) of weak Sobolev spaces it comes 
$$\qquad \qquad \|\mathcal{J}^{s/2}f \|^{q}_{L^{q, \infty}(\omega)}\leq  C \|\nabla f \|_{L^{1}(\omega)}$$ 
which is the desired result.\begin{flushright}{$\blacksquare$}\end{flushright} 
\subsection{Strong inequalities}
When $s=0$ in the weak inequalities it is possible to obtain stronger estimations. To achieve this, we will need an intermediate step:  
\begin{Proposition}\label{escalera2} 
Let $1<q<+\infty$, $\theta=\frac{1}{q}$ and $\beta=\theta/(1-\theta)$.  Then we have
$$\|f\|_ { L^{q}(\omega)}\leq C \|\nabla f\|_{L^{1}(\omega)}^{\theta } \|f\|^{1-\theta}_{\dot{B}^{-\beta, \infty}_{\infty}}$$ 
when the three norms in this inequality are bounded.  
\end{Proposition} 
\textit{\textbf{Proof}.}  We will follow closely \cite{Ledoux}. Just as in the preceding theorem, we will start by supposing that $\|f\|_{\dot{B}^{-\beta, \infty}_{\infty}}\leq 1$.  Thus, we must show the estimate 
\begin{equation}\label{ForteG} 
\|f\|_{L^{q}(\omega)}\leq C \|\nabla f\|_{L^{1}(\omega)}^{\theta}. 
\end{equation} 
Let us fix $t$ in the following way:  $t_{\alpha}=\alpha^{-2(q-1)/q}$ where $\alpha>0$.  We have then, by the thermic definition of Besov spaces, the estimate $\|H_{t}f \|_{L^\infty}\leq \alpha$. We use now the characterization of Lebesgue space given by the distribution function:  
\begin{equation}\label{formulaG} 
\frac{1}{5^{q }} \|f\|_{L^{q}(\omega)}^{q}=\int_{0}^{+\infty}\omega\left(\{x\in \mathbb{G}: |f(x)|> 5\alpha\}\right)d(\alpha^{q}).  
\end{equation} 
It now remains to estimate $\omega(\{x\in \mathbb{G}: |f(x)|> 5\alpha\})$ and for this we introduce the following thresholding function:  
\begin{equation*}\label{teta}
\Theta_{\alpha}(t)=\left\lbrace\begin{array}{l} \Theta_{\alpha}(-t)=-\Theta_{\alpha}(t)\\[4mm] 
0 \qquad \qquad\qquad\mbox{ if }\qquad 0\leq T \leq \alpha\\[4mm]
t-\alpha \qquad\qquad\mbox{ if } \qquad\alpha\leq T \leq M\alpha \\[4mm]
(M-1)\alpha \qquad\mbox{ if } \qquad T > M\alpha 
\end{array}\right.\\[3mm ] 
\end{equation*}
Here, $M$ is a parameter which depends on $q$ and which we will suppose for the moment larger than 10.  \\

This cut-off function enables us to define a new function $f_{\alpha}=\Theta_{\alpha}(f)$. We write in the next lemma some significant properties of this function $f_{\alpha}$: 
\begin{Lemme}\label{Taboo} 
\begin{enumerate} 
\item[]
\item[1)] the set defined by $\{x\in \mathbb{G}: |f(x)|> 5\alpha\}$ is included in the set $\{x\in \mathbb{G}: |f_{\alpha }(x)|> 4\alpha\}$. \\  
\item[2)] On the set $\{x\in \mathbb{G}: |f(x)|\leq M\alpha\}$ one has the estimate $|f-f_{\alpha }|\leq \alpha$. \\  
\item[3)] If $f\in \mathcal{C}^{1}(\mathbb{G})$, one has the equality $\nabla f_{\alpha}=(\nabla f)\mathds{1}_{\{\alpha\leq|f|\leq M\alpha\}}$ 
almost everywhere.  
\end{enumerate} 
\end{Lemme} 
For a proof see \cite{Ledoux}. \\  

Let us return now to (\ref{formulaG}). By the first point of the lemma above we have
\begin{equation}\label{cray} 
\int_{0}^{+\infty}\omega\left(\left\{x\in \mathbb{G}: |f(x)|> 5\alpha\right\}\right)d(\alpha^{q})\leq \int_{0}^{+\infty } \omega\left(\{x\in \mathbb{G}: |f_{\alpha }(x)|> 4\alpha\}\right)d(\alpha^{q})=I. 
\end{equation} 
We note $A_{\alpha}=\{x\in \mathbb{G}: |f_{\alpha }(x)|> 4\alpha\}$, $B_{\alpha}=\{x\in \mathbb{G}: |f_{\alpha}(x)-H_{t_{\alpha}}(f_{\alpha})(x)|> \alpha\}$ and $C_{\alpha}=\{x\in \mathbb{G}: |H_{t_{\alpha}}(f_{\alpha}-f)(x)|> 2\alpha\}$. Now, by linearity of $H_{t}$ we can write: $f_{\alpha}=f_{\alpha}-H_{t_{\alpha}}(f_{\alpha})+H_{t_{\alpha}}(f_{\alpha}-f)+H_{t_{\alpha}}(f)$. Then, holding in account the fact $\|H_{t}f\|_{L^\infty}\leq \alpha$, we obtain $A_{\alpha}\subset B_{\alpha}\cup C_{\alpha}$. Returning to (\ref{cray}), this set inclusion gives us the following inequality 
\begin{equation}\label{dosG} 
I\leq \int_{0}^{+\infty}\omega\left(B_{\alpha}\right)d(\alpha^{q})+\int_{0}^{+\infty}\omega\left(C_{ \alpha}\right)d(\alpha^{q }) 
\end{equation} 
We will study and estimate these two integrals, which we will call $I_{1}$ and $I_{2}$ respectively, by the two following lemmas:  
\begin{Lemme} For the first integral of (\ref{dosG}) we have the estimate:  
\begin{equation}\label{ARG} 
I_{1 } = \int_{0}^{+\infty}\omega\left(B_{\alpha}\right)d(\alpha^{q})\leq C\, q\log(M) \|\nabla f\|_ { L^{1}(\omega) } 
\end{equation} 
\end{Lemme} 
\textbf{\emph{Proof}.} Tchebytchev's inequality implies 
\begin{equation*} 
\omega\left(B_{\alpha}\right)\leq \alpha^{-1}\int_{\mathbb{G } }|f_{\alpha}(x)-H_{t_{\alpha}}(f_{\alpha })(x)|\omega(x)dx.  
\end{equation*} 
Using the modified Poincar\'e pseudo-inequality (\ref{CHAME1}) with $s=0$ in the above integral we obtain:  
$$\omega\left(B_{\alpha}\right)\leq C \, \alpha^{-1 } \, t_{\alpha}^{1/2}\int_{\mathbb{G } }|\nabla f_{\alpha }(x)|\omega(x)dx$$ 
Remark that the choice of $t_{\alpha}$ fixed before gives $t_{\alpha}^{1/2}=\alpha^{1-q}$, then we have
$$\omega\left(B_{\alpha}\right)\leq C \, \alpha^{-q}\int_{\{\alpha\leq|f|\leq M\alpha \} }|\nabla f(x)|\omega(x)dx.$$ 
We integrate now the preceding expression with respect to $d(\alpha^{q})$:  
$$I_{1}\leq C\int_{0}^{+\infty}\alpha^{-q}\left(\int_{\{\alpha\leq|f|\leq M\alpha \} }|\nabla f(x)|\omega(x)dx\right)d(\alpha^{q }) = C\;q\int_{\mathbb{G } }|\nabla f(x)|\left(\int_{\frac{|f|}{M}}^{|f|}\frac{d\alpha}{\alpha}\right)\omega(x)dx$$ 
It follows then $I_{1}\leq C\, q\, \log(M) \|\nabla f \|_ { L^{1}(\omega)}$ and one obtains the estimation needed for the first integral.
\begin{flushright}{$\blacksquare$}\end{flushright} 
\begin{Lemme}\label{papas}
For the second integral of (\ref{dosG}) one has the following result:  
\begin{equation*}
I_{2}=\int_{0}^{+\infty}\omega(C_{\alpha})d(\alpha^{q})\leq \frac{q}{q-1}\;\frac{1}{M^{q-1 } } \|f\|_ {L^q}^{q } 
\end{equation*} 
\end{Lemme} 
\textbf{\emph{Proof}.} For the proof of this lemma, we write:  $$|f-f_{\alpha }|=|f-f_{\alpha }|\mathds{1}_{\{|f|\leq M \alpha\}}+|f-f_{\alpha }|\mathds{1}_{\{|f|> M\alpha\}}.$$ 
As the distance between $f$ and $f_{\alpha}$ is lower than $\alpha$ on the set $ \{x\in \mathbb{G}: |f(x)|\leq M \alpha\}$, one has the inequality 
$$|f-f_{\alpha }|\leq\alpha+|f|\mathds{1}_{\{|f|> M \alpha\}}$$ 
By applying the heat semi-group to both sides of this inequality we obtain $H_{t_{\alpha}}(|f-f_{\alpha }|)\leq \alpha + H_{t_{\alpha}}(|f|\mathds{1}_{\{|f|> M \alpha\}})$ and we have then the following set inclusion $C_{\alpha}\subset \left\{x\in \mathbb{G}: H_{t_{\alpha}}(|f|\mathds{1}_{\{|f|> M \alpha\}})>\alpha\right\}$.
Thus, considering the measure of these sets and integrating with respect to $d(\alpha^{q})$, it comes 
\begin{equation*} 
I_{2}=\int_{0}^{+\infty}\omega\left(C_{\alpha}\right)d(\alpha^{q})\leq\int_{0}^{+\infty}\omega\bigg(\{H_{t_{\alpha}}(|f|\mathds{1}_{\{|f|> M\alpha\}})>\alpha\}\bigg)d(\alpha^{q}) 
\end{equation*} 
We obtain now, by applying Tchebytchev inequality, the estimate 
$$I_{2}\leq\int_{0}^{+\infty}\alpha^{-1}\bigg(\int_{\mathbb{G}}H_{t_{\alpha}}\left(|f|\mathds{1}_{\{|f|> M\alpha\}}\right)\omega(x)dx\bigg)d(\alpha^{q}),$$ 
then by Fubini's theorem we have
$$I_{2}\leq q\int_{\mathbb{G}}|f(x)|\bigg(\int_{0}^{+\infty}\mathds{1}_{\{|f|> M\alpha\}}\alpha^{q-2}d\alpha\bigg)\omega(x)dx=\frac{q}{q-1}\int_{\mathbb{G}}|f(x)|\frac{|f(x)|^{q-1}}{M^{q-1}}\omega(x)dx= \frac{q}{q-1}\frac{1}{M^{q-1 }}\|f\|_{L^{q}(\omega)}^{q}.$$ 
And this concludes the proof of this lemma.
\begin{flushright}{$\blacksquare$}\end{flushright} 
We finish the proof of proposition \ref{escalera2} by connecting together these two lemmas \textit{i.e.}:
$$\frac{1}{5^{q}} \|f\|_ {L^q(\omega)}^{q}\leq Cq\, \log(M) \|\nabla f\|_ {L^1(\omega)}+\frac{q}{q-1}\frac{1}{M^{q-1 } } \|f\|_{L^q(\omega)}^{q}$$ 
Since we supposed all the norms bounded and  $M\gg 1$, we finally have
$$\left(\frac{1}{5^{q}}-\frac{q}{q-1}\frac{1}{M^{q-1}}\right) \|f\|_ {L^q(\omega)}^{q}\leq C q\, \log(M) \|\nabla f\|_{L^1(\omega)}$$ 
\begin{flushright}{$\blacksquare$}\end{flushright} 
The proof of the theorem \ref{smile00} is not yet completely finished. The last step is provided by the 
\begin{Proposition}\label{escalera3} 
In proposition \ref{escalera2} it is possible to consider only the two assumptions 
$\nabla f\in L^{1}(\mathbb{G},\omega)$ and $f\in \dot{B}^{-\beta, \infty}_{\infty}(\mathbb{G})$.  
\end{Proposition}
\textbf{\emph{Proof}.} For the proof of this proposition we will build a homogeneous Littlewood-Paley like approximation of $f$ writing:  
\begin{equation*} 
f_{j}=\left(\int_{0}^{+\infty}\left(\varphi(2^{-2j}\lambda)-\varphi(2^{2j}\lambda)\right)dE_{\lambda}\right)(f) \qquad (j\in \mathbb{N})
\end{equation*} 
where $\varphi$ is a $\mathcal{C}^{\infty}(\mathbb{R}^+)$ function such that $\varphi=1$ on $]0,1/4[$ and $\varphi=0$ on $[1, +\infty[$.
\begin{Lemme} 
If $q>1$, if $\nabla f\in L^{1}(\mathbb{G}, \omega)$ and if $f\in \dot{B}^{-\beta, \infty}_{\infty}(\mathbb{G})$ then $\nabla f_j\in L^{1}(\mathbb{G}, \omega)$, $f_j\in \dot{B}^{-\beta, \infty}_{\infty}(\mathbb{G})$ and $f_{j}\in L^{q}(\mathbb{G},\omega)$.  
\end{Lemme} 
\textbf{\emph{Proof}.} The fact that $\nabla f_j\in L^{1}(\mathbb{G}, \omega)$ and $f_j\in \dot{B}^{-\beta, \infty}_{\infty}(\mathbb{G})$ is an easy consequence of the definition of $f_j$. For $f_{j}\in L^{q}(\mathbb{G},\omega)$ the starting point is given by the relation:  
\begin{equation*}
f_{j}=\left(\int_{0}^{+\infty}m(2^{-2j}\lambda)\, dE_{\lambda}\right)2^{-2j}\mathcal{J}(f), 
\end{equation*} 
where we noted $$m(2^{-2j}\lambda)=\frac{\varphi(2^{-2j}\lambda)-\varphi(2^{2j}\lambda)}{2^{-2j}\lambda}.$$ 
Observe that the function $m$ vanishes near of the origin and satisfies the assumptions of corollary \ref{black}. We obtain then the following identity where $M_{j}\in \mathcal{S}(\mathbb{G})$ is the kernel of the operator $m(2^{-2j}\mathcal{J})$:  
$$f_{j}=2^{-2j}\mathcal{J}f\ast M_{j}=2^{-2j}\nabla f\ast \tilde{\nabla}M_{j},$$ 
where we denoted $\tilde{\nabla}$ the gradient of the right invariant vectors fields and used the property (\ref{izder}). Let us now calculate the norm $L^{q}(\mathbb{G},\omega)$ in the preceding identity:
$$\|f_{j}\|_{L^{q}(\omega)}=\|2^{-2j}\nabla f\ast\tilde{\nabla}M_{j}\|_{L^{q}(\omega)}\leq 2^{-2j}\|\nabla f\|_{L^{1}(\omega)}\|\tilde{\nabla}M_{j}\|_{L^{q}(\omega)}.$$ 
Finally, we obtain:  
$$\|f_{j } \|_ { L^{q}(\omega)}\leq C\, 2^{j(N(1-\frac{1}{q})-1) } \|\nabla f \|_ { L^{1}(\omega)}<+\infty$$ 
\begin{flushright}{$\blacksquare$}\end{flushright} 

Thanks to this estimate, we can apply the proposition \ref{escalera2} to $f_{j}$ whose $L^{q}(\mathbb{G}, \omega)$ norm is bounded, and we obtain:  
$$\|f_{j } \|_ { L^{q}(\omega)}\leq C \|\nabla f_j\|_ { L^{1}(\omega)}^{\theta } \|f_j\|_ { \dot{B}^{-\beta, \infty}_{\infty}}^{1-\theta}.$$ 
Now, since $f\in \dot{B}^{-\beta, \infty}_{\infty}(\mathbb{G})$, we have $f_{j } \rightharpoonup f$ in the sense of distributions. It follows  
$$\|f\|_{ L^{q}(\omega)}\leq \underset{j \to +\infty}{\lim \inf } \|f_{j } \|_ { L^{q}(\omega)}\leq C \|\nabla f\|_ { L^{1}(\omega)}^{\theta }
 \|f\|_ { \dot{B}^{-\beta, \infty}_{\infty}}^{1-\theta}.$$ 
We restricted ourselves to the two initial assumptions, namely $\nabla f\in L^{1}(\mathbb{G}, \omega)$ and $f\in\dot{B}^{-\beta, \infty}_{\infty}(\mathbb{G})$. The strong inequalities (\ref{smile0}) are now completely proved for stratified groups.  
\begin{flushright}{$\blacksquare$}\end{flushright} 

\subsection{Maximal function and Improved Sobolev inequalities}
We will study now theorem \ref{PGLR}. Just as for weak inequalities (\ref{ForIsoLap}), we can rewrite (\ref{PGLR1}) in the following way 
\begin{equation*}
\| \mathcal{J}^{\frac{s-s_1}{2}}f\|_{L^q(\omega)}\leq C \|f\|_{L^p(\omega)}^\theta\|f\|_{\dot{B}^{-\beta-s_1,\infty}_{\infty}}^{1-\theta}
\end{equation*}
where $1<p<q<+\infty$, $\theta=p/q$, $s=\theta s_1 -(1-\theta)\beta$ and $-\beta<s<s_1$. Using the sub-Laplacian fractional powers characterization  we have the identity
\begin{equation}\label{Kashmor}
\mathcal{J}^{\frac{-\alpha}{2}}f(x)=\frac{1}{\Gamma(\frac{\alpha}{2})}\int_{0}^{+\infty}t^{\frac{\alpha}{2}-1}H_t f(x)dt=\frac{1}{\Gamma(\frac{\alpha}{2})}\left(\int_{0}^{T}t^{\frac{\alpha}{2}-1}H_t f(x)dt+\int_{T}^{+\infty}t^{\frac{\alpha}{2}-1}H_t f(x)dt\right)
\end{equation}
where $\alpha=s_1-s>0$ and $T$ will be fixed in the sequel.\\

For studying each one of these integrals we will use the estimates
\begin{enumerate}
\item[$\bullet$] $|H_tf(x)|\leq C \mathcal{M}_B f(x)$ \qquad\qquad\qquad\;(by lemma \ref{Led11})\\

\item[$\bullet$] $|H_tf(x)|\leq C t^{\frac{-\beta-s_1}{2}}\|f\|_{\dot{B}^{-\beta-s_1, \infty}_{\infty}}$ \qquad (by the thermic definition of Besov spaces (\ref{BesovChaleurGhomo}))\\
\end{enumerate}
Then, applying these inequalities in (\ref{Kashmor}) we obtain
$$|\mathcal{J}^{\frac{-\alpha}{2}}f(x)|\leq \frac{c_1}{\Gamma(\frac{\alpha}{2})}T^{\frac{\alpha}{2}} \mathcal{M}_B f(x)+ \frac{c_2}{\Gamma(\frac{\alpha}{2})}T^{\frac{\alpha-\beta-s}{2}}\|f\|_{\dot{B}^{-\beta-s_1, \infty}_{\infty}}.$$
We fix now $$T=\left(\frac{\|f\|_{\dot{B}^{-\beta-s_1, \infty}_{\infty}}}{\mathcal{M}_B f(x)}\right)^{ \frac{2}{\beta+s_1}}$$
and we get 
$$|\mathcal{J}^{\frac{-\alpha}{2}}f(x)|\leq \frac{c_1}{\Gamma(\frac{\alpha}{2})}\mathcal{M}_B f(x)^{1-\frac{\alpha}{\beta+s_1}}+ \frac{c_2}{\Gamma(\frac{\alpha}{2})}\mathcal{M}_B f(x)^{1-\frac{\alpha}{\beta+s_1}}\|f\|^{\frac{\alpha}{\beta+s_1}}_{\dot{B}^{-\beta-s_1, \infty}_{\infty}}.$$
Since $\frac{\alpha}{\beta+s_1}=1-\theta$, we have
$$|\mathcal{J}^{\frac{-\alpha}{2}}f(x)|\leq \frac{c}{\Gamma(\frac{\alpha}{2})}\mathcal{M}_B f(x)^{\theta}\|f\|^{1-\theta}_{\dot{B}^{-\beta-s_1, \infty}_{\infty}}.$$
Multiplying this inequality by a $A_p$ weight $\omega$, using the fact that $A_p\subset A_q$ if $p<q$, and the property of maximal function (\ref{MaximalProper}) we obtain
$$\|\mathcal{J}^{\frac{-\alpha}{2}}f(x)\|_{L^q(\omega)}\leq c\|f\|_{L^p(\omega)}^{\theta}\|f\|^{1-\theta}_{\dot{B}^{-\beta-s_1, \infty}_{\infty}}$$
and we are done. 
\begin{flushright}{$\blacksquare$}\end{flushright} 
It is worth noting the simplicity of this proof: the arguments used are classical tools from harmonic analysis. 

\begin{flushright}
\begin{minipage}[r]{80mm}
Diego \textsc{Chamorro}\\[5mm]
Laboratoire d'Analyse et de Probabilités\\ 
ENSIIE\\[2mm]
1 square de la résistance,\\
91025 Evry Cedex\\[2mm]
diego.chamorro@m4x.org
\end{minipage}
\end{flushright}


\begin{thebibliography}{2}
\bibitem{Bahouri}
H. \textsc{Bahouri}, P. \textsc{G\'erard} \& Ch. \textsc{XU}.
\emph{Espaces de Besov et estimations de Strichartz g\'en\'eralis\'ees sur le groupe de Heisenberg.}
Journal d'Analyse Math\'ematique, Vol. 82 (2000).
\bibitem{Bergh}
J. \textsc{Bergh} \& J. \textsc{L\"ofst\"rom}. \emph{Interpolation Spaces}. 
Grundlehren der mathematischen Wissenschaften, 223. Springer Verlag (1976).
\bibitem{Chame}
D. \textsc{Chamorro}. \textit{In\'egalit\'es de Gagliardo-Nirenberg pr\'ecis\'ees sur le groupe de Heisenberg}, Th\`ese de troisi\`eme cycle.  ENS Cachan (2006).
\bibitem{Cohen2}
A. \textsc{Cohen}, W. \textsc{Dahmen}, I. \textsc{Daubechies} \& R. \textsc{De Vore}. \emph{Harmonic Analysis of the space BV}. Rev. Mat. Iberoamericana 19, n$^\circ$1, 235-263 (2003).
\bibitem{Folland}
G. \textsc{Folland}. \emph{Subelliptic estimates and function spaces on nilpotent Lie groups}. Ark. Mat. 13, 161-208 (1975).
\bibitem{Folland2}
G. \textsc{Folland} \& E. M. \textsc{Stein}. \emph{Hardy
Spaces on homogeneous groups}. Mathematical Notes, 28, Princeton University Press (1982).
\bibitem{Furioli2}
G. \textsc{Furioli}, C. \textsc{Melzi} \& A. \textsc{Veneruso}. \emph{Littlewood-Paley decomposition and Besov spaces on Lie groups of polynomial growth.}
Math. Nachr. 279, n$^\circ$ 9-10, 1028-1040 (2006).
\bibitem{Garcia}
J. \textsc{Garc\'ia-Cuerva} \& J.L. \textsc{Rubio de Francia}. \emph{Weighted norm inequalities and related topics}. 
Mathematics studies 116, North Holland (1985).
\bibitem{Garofalo}
N. \textsc{Garofalo} \& D \textsc{Nhieu}. \emph{Isoperimetric and Sobolev inequalities for Carnot-Carath\'eodory spaces and the existence of minimal surfaces}. Communications on Pure and Applied Mathematics, Vol XLIX, 1081-1144 (1996).
\bibitem{GMO}
P. \textsc{G\'erard}, Y. \textsc{Meyer} \& F. \textsc{Oru}. \emph{In\'egalit\'es de Sobolev Pr\'ecis\'ees}. Equations aux 
D\'eriv\'ees Partielles, S\'eminaire de l'Ecole Polytechnique, expos\'e n$^\circ$ IV (1996-1997).
\bibitem{Grafakos}
L. \textsc{Grafakos}. \emph{Classical and Modern Fourier Analysis}. Prentice Hall (2004).
\bibitem{Hulanicki}
A. \textsc{Hulanicki}. \emph{A functional calculus for Rockland operators on nilpotent Lie groups}. Studia Mathematica T. LXXVIII (1984).
\bibitem{Ledoux}
M. \textsc{Ledoux}. \emph{On improved Sobolev embedding theorems}. Math. Res. Letters 10, 659-669 (2003).
\bibitem{Meyer}
Y. \textsc{Meyer}. \emph{Ondelettes et Op\'erateurs}. Hermann (1990).
\bibitem{Saka}
K. \textsc{Saka}. \emph{Besov Spaces and Sobolev spaces on a nilpotent Lie group}. 
Thoku. Math. Journ.  Vol. 31, p. 383-437 (1979).
\bibitem{Stein0}
E. M. \textsc{Stein}. \emph{Topics in Harmonic analysis}. Annals of mathematics studies, 63. Princeton University Press (1970).
\bibitem{Stein2}
E. M. \textsc{Stein}. \emph{Harmonic Analysis}. Princeton University Press (1993).
\bibitem{Strichartz}
R. \textsc{Strichartz}. \emph{Self-similarity on nilpotent Lie groups}. Contemporary Mathematics, volume 142, 123-157 (1992).
\bibitem{Triebel}
H. \textsc{Triebel}. \emph{Theory of function spaces II}.
Birkh\"auser (1992).
\bibitem{Varopoulos}
N. Th. \textsc{Varopoulos}, L. \textsc{Saloff-Coste} \& T. \textsc{Coulhon}. \emph{Analysis and geometry on groups}.
Cambridge Tracts in Mathematics, 100 (1992).\\[10mm]
\end{thebibliography}
\end{document}